\documentclass[12pt]{article}

\usepackage[dvips]{graphicx}
\usepackage[utf8]{inputenc}
\usepackage{psfrag}
\usepackage{amsthm}
\usepackage{amsmath}
\usepackage{amssymb}

\usepackage[agsm]{harvard}

\newtheorem{Theorem}{Theorem}
\newtheorem{Definition}{Definition}
\newtheorem{Lemma}{Lemma}

\newcommand{\abs}[1]{\left\vert #1 \right\vert}
\newcommand{\norm}[1]{\left\Vert #1 \right\Vert}

\newcommand{\R}{{\mathbb R}}  
\newcommand{\E}{{\mathbb E}}  

\citationmode{abbr}
 
\title{Global convergence of quorum-sensing networks} \author{Giovanni
  Russo\thanks{G. Russo is with the Department of Systems and Computer
    Engineering, University of Naples Federico II, Italy. Work done
    while visiting the Nonlinear Systems Laboratory, Massachusetts
    Institute of Technology. E-mail: {\tt\small
      giovanni.russo2@unina.it}} \ \ and \ \ Jean-Jacques
  Slotine\thanks{J.J. Slotine is with the Nonlinear Systems
    Laboratory, Massachusetts Institute of Technology, United
    States. E-mail: {\tt\small jjs@mit.edu}}}

\begin{document}

\markboth{Convergence of quorum-sensing  networks}{Convergence of quorum-sensing  networks}

\maketitle

\begin{abstract}

In many natural synchronization phenomena, communication between
individual elements occurs not directly, but rather through the
environment. One of these instances is bacterial quorum sensing, where
bacteria release signaling molecules in the environment which in turn
are sensed and used for population coordination.  Extending this
motivation to a general nonlinear dynamical system context, this paper
analyzes synchronization phenomena in networks where communication and
coupling between nodes are mediated by shared dynamical quantities,
typically provided by the nodes' environment. Our model includes the
case when the dynamics of the shared variables themselves cannot be
neglected or indeed play a central part. Applications to examples from
systems biology illustrate the approach.

\end{abstract}


\section{Introduction}

Many dynamical phenomena in biology involve some form of {\it
  synchronization}.  Synchronization has attracted much research both
from the theoretical, see e.g. \cite{Str_03},\cite{You_Cox_Wei_Arn_04}, \cite{McM_Kop_Has_Col_02} to cite just a few, and experimental
\cite{Yag_Ise_Mat_Oku_Yag_03}, \cite{Pye_69} viewpoints.  The particular case
of synchronized {\it time-periodic} processes, where time-scales can
range from a few milliseconds to several years
\cite{Win_01,New_Bar_Wat_06}, includes e.g.  circadian rhythms in
mammals \cite{Gon_Ber_Wal_Kra_Her_05}, the cell cycle
\cite{Tys_Csi_Now_02}, spiking neurons \cite{Izh_06} and
respiratory oscillations \cite{Hen_04}.


When modelling such networks, it is often assumed that each node
communicates directly with other nodes in the network, see
e.g. \cite{Par_Fen_Dur_08,Boh_Oja_08} and references therein.  In many
natural instances, however, network nodes do not communicate directly,
but rather by means of noisy and continuously changing environments.
Bacteria, for instance, produce, release and sense signaling molecules
(so-called autoinducers) which can diffuse in the environment and are
used for population coordination. This mechanism, known as {\it quorum
  sensing}~\cite{Mil_Bas_01,Nar_Bas_Lev_08,Ng_Bas_09} is believed to
play a key role in bacterial infection, as well as e.g. in
bioluminescence and biofilm formation \cite{Ane_Pir_Jun_09},
\cite{Nad_Xav_Lev_Fos_08}.  In a neuronal context, a mechanism similar
to that of quorum sensing may involve \emph{local field potentials},
which may play an important role in the synchronization of clusters of
neurons,
\cite{Per_Pez_Sah_Mit_And_02,Fregnac_09,Tab_Slo_Pha_09,Ana_Mon_Bar_Buz_Koc_10}. 

From a network dynamics viewpoint, the key
characteristic of quorum sensing-like mechanisms lies in the fact
that communication between nodes (e.g. bacteria) occurs by means of a
shared quantity (e.g. autoinducer concentration). Furthermore, the
production and degradation rates of such a quantity are affected by
all the nodes of the network. Therefore, a detailed model of such a
mechanism needs to keep track of the temporal evolution of the shared
quantity, resulting in an additional set of ordinary differential
equations.

Mathematical work on such quorum sensing topologies is relatively
sparse (e.g., \cite{Gar_Elo_Str_04,Tab_Slo_Pha_09,Rus_diB_09b,Kat_08}) compared to that
on diffusive topologies, and often neglects quorum variable dynamics
or the dynamics of the environment. This sparsity of results is
somewhat surprising given that, besides its biological pervasiveness,
quorum sensing may also be viewed as an astute ``computational''
tool. Specifically, use of a shared variable in effect significantly
reduces the number of links required to achieve a given
connectivity~\cite{Tab_Slo_Pha_09}.

In this paper, we derive sufficient conditions for the coordination of
nodes communicating through dynamical quorum sensing mechanisms. These
results can be used both to analyse natural networks, and to guide
design of communication mechanisms in synthetic or partially synthetic
networks.  We first consider, in Section \ref{sec:basic_model}, the
case where the network nodes (e.g., the biological entities populating
the environment) are all identical or nearly identical. We then focus,
in Section \ref{sec:basic_model_multiple}, on networks composed of
heterogeneous nodes, i.e., nodes of possibly diverse dynamics. In this
case we provide sufficient conditions ensuring that all the network
nodes sharing the same dynamics converge to a common behavior, a
particular instance of so-called concurrent synchronization
\cite{Pha_Slo_07}.  In Section~\ref{sec:cluster}, the results are
further extended to a distributed version of quorum sensing, where
multiple groups of possibly heterogeneous nodes communicate by means
of multiple media.  Finally, in Section~\ref{sec:synchro_control}, we
propose a strategy for controlling the common asymptotic evolution of
the network nodes. Section~\ref{N} studies the dependence of
synchronization properties on the number of nodes, a question of
interest e.g. in the context of cell
proliferation. Section~\ref{sec:examples} illustrates the general
approach with a set of examples.

Our proofs are based on nonlinear contraction theory
(\cite{Loh_Slo_98}), a viewpoint on incremental stability which we
briefly review in Section~\ref{sec:math_tools}, and which has emerged
as a powerful tool in applications ranging from Lagrangian mechanics
to network control. Historically, ideas closely related to contraction
can be traced back to~\cite{Hartmann} and even to~\cite{Lewis} (see
also~\cite{Pav_Pog_Wou_Nij,Ang_02}, and e.g.~\cite{pde} for a more
exhaustive list of related references).  As pointed out
in~\cite{Loh_Slo_98}, contraction is preserved through a large variety
of systems combinations, and in particular it represents a natural
tool for the study and design of nonlinear state observers, and by
extension, of synchronization mechanisms~\cite{Wan_Slo_05}.


\section{Contraction theory tools}\label{sec:math_tools}

\subsection{Basic results}\label{sec:basic_res}

The basic result of nonlinear contraction analysis \cite{Loh_Slo_98}
which we shall use in this paper can be stated as follows.

\begin{Theorem}[Contraction]
\label{theorem:contraction}
Consider the $m$-dimensional deterministic system
\begin{equation}
\label{equ:main}
\dot x = f(x,t)
\end{equation}
where $f$ is a smooth nonlinear function.  The system is said to be
{\it contracting} if any two trajectories, starting from different
initial conditions, converge exponentially to each other. A sufficient
condition for a system to be contracting is that there exists a
constant invertible matrix $\Theta$ such that the so-called
generalized Jacobian
\begin{equation}
\label{equ:generalized_Jacobian}
F(x,t)  =  \Theta \
\frac{\partial f}{\partial x} (x,t) \ \Theta^{-1}
\end{equation}
verifies
$$
\exists \lambda > 0, \ \forall x, \ \forall t \ge 0, \ \ \mu(F(x,t)) \ \le \ - \lambda
$$
where $\mu$ is one the the standard matrix measures in Table 1.
The scalar $\lambda$ defines the contraction rate of the system.
\end{Theorem}

For convenience, in this paper we will also say that a {\it function}
$f(x,t)$ is contracting if the system $\ \dot x= f(x,t) \ $
satisfies the sufficient condition above.  Similarly, we
will then say that the corresponding Jacobian {\it matrix} $\ \frac{\partial
f}{\partial x}(x,t)$ is contracting.

\begin{table}[th] 
\caption{Standard Matrix measures}
\centering
\label{tab:matrix_measures}
\begin{tabular}{|c| c|} 
\hline
vector norm, $\abs{\cdot}$ & induced matrix measure, $\mu\left(A\right)$\\
\hline
$\abs{x}_1= \sum_{j=1}^m\abs{x_j}$ & $ \mu_1\left(A\right)= \max_{j} \left( a_{jj}+\sum_{i \ne j}  \abs{a_{ij}} \right)$\\
\hline
$\abs{x}_2= \left( \sum_{j=1}^n\abs{x_j}^2\right)^{\frac{1}{2}}$ & $ \mu_2\left(A\right)=\max_{i} \left( \lambda_i\left\{\frac{A+A^\ast}{2}\right\}\right)$\\
\hline
$\abs{x}_\infty= \max_{1 \le j \le m} \abs{x_j}$ & $\mu_{\infty}\left(A\right)= \max_{i} \left( a_{ii}+\sum_{j \ne i} \mid a_{ij}\mid \right)$\\
\hline
\end{tabular}
\end{table}

We shall also use the following two properties of contracting systems,
whose proofs can be found in~\cite{Loh_Slo_98,Slo_03}.

\noindent {\bf Hierarchies of contracting systems}\ \ \ Assume that 
the Jacobian of (\ref{equ:main}) is in the form
\begin{equation}\label{eqn:hierarchy_general}
\frac{\partial f}{\partial x}(x,t)\ = \
\left[\begin{array}{*{20}c}
J_{11} & J_{12}\\
0 & J_{22}\\
\end{array}\right]
\end{equation}
corresponding to a hierarchical dynamic structure. The $J_{ii}$ may be
of different dimensions. Then, a sufficient condition for the system
to be contracting is that (i) the Jacobians $J_{11}$, $J_{22}$ are
contracting (possibly with different $\Theta$'s and for different
matrix measures), and (ii) the matrix $J_{12}$ is bounded.

\noindent {\bf Periodic inputs}\ \ \ Consider the system
\begin{equation}\label{periodic}
\dot x =f\left(x,r(t)\right)
\end{equation}
where the input vector $r(t)$ is periodic, of period $T$. Assume that
the system is contracting (i.e., that the Jacobian matrix $\
\frac{\partial f}{\partial x}(x,r(t))$ is contracting for any $r(t)$).  Then the
system state $x(t)$ tends exponentially towards a periodic state of
period $T$.

\subsection{Partial Contraction}

A simple yet powerful extension to nonlinear contraction theory is the concept
of \emph{partial} contraction~\cite{Wan_Slo_05}. 

\begin{Theorem}[Partial contraction]
\label{theorem:partial_contraction}
Consider a smooth nonlinear $m$-dimensional system of the form $\dot x= f(x,x,t)$ and assume that the so-called system $\dot y= f(y,x,t)$ is contracting with respect to $y$.  If a particular solution of
the auxiliary $y$-system verifies a smooth specific property,
then all trajectories of the original $x$-system verify this
property exponentially. The original system is said to be \emph{partially
contracting}.
\end{Theorem}

Indeed, the virtual $y$-system has two particular solutions, namely
$y(t) = x\left(t\right)$ for all $t \ge 0$ and the
particular solution with the specific property. Since all trajectories
of the $y$-system converge exponentially to a single trajectory, this
implies that $x\left(t\right)$ verifies the specific property
exponentially.

\subsection{Networks of contracting nodes}\label{sec:preliminaries}

This section introduces preliminary results on concurrent synchronization of networks,
which will be used in the rest of the paper. 

Consider a network consisting of $N$ heterogeneous nodes:
\begin{equation}\label{eqn:net_multiple_dynamics}
\dot x_i =f_{\gamma\left(i\right)}\left(x_i,t\right) + \sum_{j\in N_i}\left[h_{\gamma\left(i\right)}\left(x_j\right)-h_{\gamma\left(i\right)}\left(x_i\right)\right]
\end{equation}
where $N_i$ denotes the set of neighbors of node $i$ and $\gamma$ is a function defined between two set of indices (not necessarily a permutation), i.e.
\begin{equation}\label{eqn:gamma}
\gamma : \left\{1,\ldots,N\right\} \rightarrow \left\{1,\ldots, s\right\} \quad s\le N
\end{equation}
 Thus, two nodes of (\ref{eqn:net_multiple_dynamics}), e.g. $x_i$ and
 $x_j$, share the same dynamics and belong to the $p$-th group (denoted with $\mathcal{G}_p$), i.e. $x_i, x_j \in \mathcal{G}_p$, if and
 only if $\gamma\left(i\right)=\gamma\left(j\right) = p$. The dimension of the nodes' state variables belonging to group $p$ is $n_{\gamma(i)}$, i.e. $x_i\in \R^{n_{\gamma(i)}}$ for any $x_i \in \mathcal{G}_p$. In what follows we assume that the Jacobian of the coupling functions $h_{\gamma(i)}$ are diagonal matrices with nonnegative diagonal elements. We will derive conditions ensuring \emph{concurrent
 synchronization} of (\ref{eqn:net_multiple_dynamics}), i.e. all nodes belonging to the same group exhibit the same regime behavior.
 
In what follows the following standard assumption (see \cite{Pha_Slo_07} and references therein) is made on the interconnections between the agents belonging to different groups, \cite{Gol_Ste_Tor_05}.

\begin{Definition}\label{ass:input_symmetry}
Let $i$ and $j$ be two nodes of a group $G_p$, and if they receive their input from elements $i'$, $j'$ respectively, then: (ii) $i'$ and $j'$ belong to the same group $G_{p'}$; (ii) the coupling functions between $i$-$i'$ and $j$-$j'$ are the same; (iii) the inputs to $i$ and $j$ coming from different groups are the same.
If these assumptions are satisfied, then nodes $i$ and $j$ are said to be \emph{input-equivalent}.
\end{Definition}

Given this definition, we can state the following theorem, which
generalizes results in \cite{Pha_Slo_07} to the case of arbitrary
norms. Its proof is provided in the Appendix.

\begin{Theorem}\label{thm:network_convergence_multiple_nodes}
Assume that in (\ref{eqn:net_multiple_dynamics}) the nodes belonging to the same group are all input-equivalent and that the nodes dynamics are all contracting. Then, all node trajectories sharing the same dynamics converge towards each other, i.e. for any $x_i$, $x_j \in \mathcal{G}_p$, $p=1,\ldots,s$,
$$
 \abs{x_j\left(t\right) - x_i\left(t\right)}\rightarrow 0\ \ {\rm as}\ \ t\rightarrow +\infty
$$
\end{Theorem}
In the case of networks of identical nodes dynamics, the above result amounts to only requiring contraction for each node.

\section{Main Results}\label{sec:main_res}

In this Section we present our main results. We first provide sufficient conditions for the synchronization of a network composed by $N$ nodes communicating over a common medium, which is characterized by some nonlinear dynamics. We then extend the analysis to a number of cases, by providing sufficient conditions for the convergence of networks composed of nodes having different dynamics (non-homogeneous nodes) and communicating over multiple (possibly non-homogeneous) media.

\subsection{The basic mathematical model and convergence analysis} \label{sec:basic_model}

In the following, we analyze the convergent behavior of the network schematically represented in Figure \ref{fig:quorum_1} (left). In such a network, the $N$ nodes are assumed to all share the same smooth dynamics and to communicate by means of the same common medium, characterized by some smooth dynamics:
\begin{equation}\label{eqn:general_model}
\begin{array}{*{20}l}
\dot x_ i = f\left(x_i,z,t\right) & i =1, \ldots, N\\
\dot z = g\left(z, \Psi\left(x_1,\ldots, x_N\right),t\right)
\end{array}
\end{equation}
A simplified version of the above model was recently analyzed by means of a graphical algorithm in \cite{Rus_diB_09c}. In the above equation, the set of state variables of the nodes is $x_i$, while the set of the state variables of the common medium
dynamics is $z$. Notice that the nodes dynamics and the medium dynamics can be of different dimensions (e.g. $x_i \in\R^n$, $z \in \R^d$). The dynamics of the nodes affect the dynamics of the common medium by means of some (coupling, or input) function, $\Psi:\R^{Nn}\rightarrow \R^d$. These  functions may 
depend only on some of the components of the $x_i$ or of $z$ (as 
the example in Section \ref{sec:synchro_example} illustrates).

The following result is a sufficient condition for convergence of all nodes trajectories of (\ref{eqn:general_model}) towards each other.

\begin{Theorem}\label{thm:quorum_convergence}
All nodes trajectories of network (\ref{eqn:general_model}) globally exponentially converge towards each other if the function $f\left(x,v(t),t\right)$ is contracting for any $v(t)\in \R^d$.
\end{Theorem}
\proof
The proof is based on partial contraction (Theorem~\ref{theorem:partial_contraction}). Consider the following \emph{reduced order} virtual system
\begin{equation}\label{eqn:virtual_quorum}
\dot y = f\left(y,z,t\right)
\end{equation}
Notice that now $z(t)$ is seen as an exogenous input to the virtual system. Furthermore, substituting $x_i$ to the virtual state variable $y$ yields the dynamics of the $i$-th node. That is, $x_i$, $i=1,\ldots,N$, are particular solutions of the virtual system. Now, if such a system is contracting, then all of its solutions will converge towards each other. Since the nodes state variables are particular solutions of (\ref{eqn:virtual_quorum}), contraction of the virtual system implies that, for any $i,j = 1,\ldots,N$:
$$
\abs{x_i-x_j}\rightarrow 0
$$
as $t\rightarrow +\infty$. 

The Theorem is proved by noting that by hypotheses the function $f(x,v(t),t)$ is contracting for any exogenous input $v(t)$. This in particular implies that $f(y,z,t)$ is contracting, i.e. (\ref{eqn:virtual_quorum}) is contracting.
\endproof

\subsection*{Remarks}
\begin{itemize}
\item The function $\Psi\left(x_1,\ldots,x_N\right)$ is often of the form
$$
\Psi\left(x_1,\ldots,x_N\right) := \sum_{i=1}^Nu(x_i)
$$
where $u:\R^n\rightarrow\R^d$ and all network nodes affect the medium
dynamics in a similar way.

\item In applications, the coupling between the nodes and the common medium is 
often assumed to be diffusive. Model (\ref{eqn:general_model}) 
then reduces to:
\begin{equation}\label{eqn:general_model_part}
\begin{array}{*{20}l}
\dot x_ i = f\left(x_i,t\right) + k_z\left( z\right) - k_x\left( x_i\right) & i =1, \ldots, N\\
\dot z = g\left(z,t\right) + \sum_{i=1}^N \left[u_x\left( x_i\right) - u_z\left( z\right)\right]
\end{array}
\end{equation}
That is, the nodes and the common medium are coupled
by means of the smooth coupling functions $k_z: \R^d\rightarrow \R^n$, $k_x: \R^n\rightarrow \R^n$ and
$u_x:\R^n\rightarrow \R^d$, $u_z:\R^d\rightarrow \R^d$. These  functions may 
depend only on some of the components of the $x_i$ or of $z$ (as 
the example in Section \ref{sec:synchro_example} illustrates). In this case, Theorem \ref{thm:quorum_convergence} implies that synchronization is attained if $f\left(x,t\right) - k_x\left( x\right)$ is contracting. Similar results are easily
derived for the generalizations of the above model presented in what follows.

\item The result also applies to the case where the quorum signal
is based not on the $x_i$'s themselves, but rather on variables deriving
from the $x_i$'s through some further nonlinear dynamics.
Consider for instance the system

\ \ \ \ \ \ \ \ \ \ ${\dot x}_i = f\left(x_i,z,t\right)  \ \ \ \ \ \ \ \ \ \ \ \ \ \ \ \ \ i =1, \ldots, N$

\ \ \ \ \ \ \ \ \ \ ${\dot r}_i = h\left(r_i,x_i,z,t\right)  \ \ \ \ \ \ \ \ \ \ \ \ \ \ i =1, \ldots, N$

\ \ \ \ \ \ \ \ \ \ ${\dot z} = g\left(z, \Psi\left(r_1,\ldots, r_N\right),t\right)$

\noindent Theorem~\ref{thm:quorum_convergence} can be applied directly
by describing each network node by the augmented state $(x_i, r_i),$
and using property~(\ref{eqn:hierarchy_general}) on hierarchical
combinations to evaluate the contraction properties of the augmented
network dynamics;
\item Similarly, each network "node" may actually be composed of several subsystems, with each subsystem synchronizing with its analogs in other nodes.
\end{itemize} 

\subsection{Multiple systems communicating over a common medium} \label{sec:basic_model_multiple}

\begin{figure}[thbp]
\begin{center}
  \includegraphics[width=7cm]{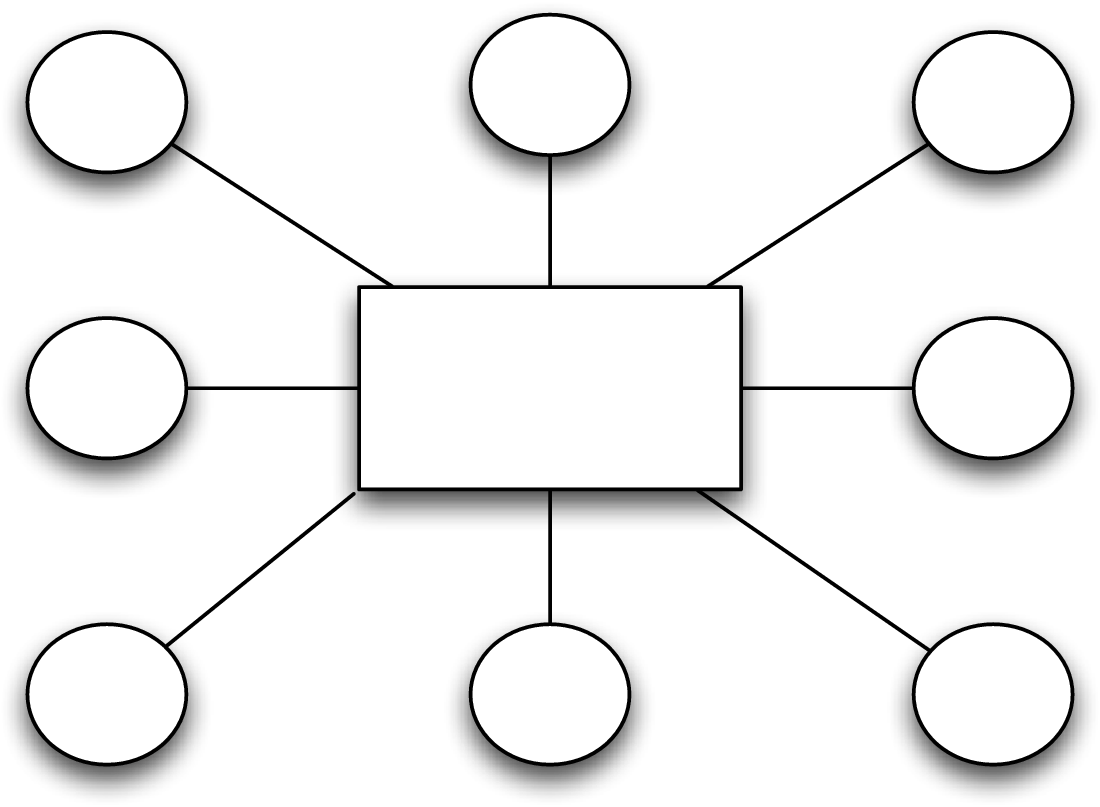}
    \includegraphics[width=7cm]{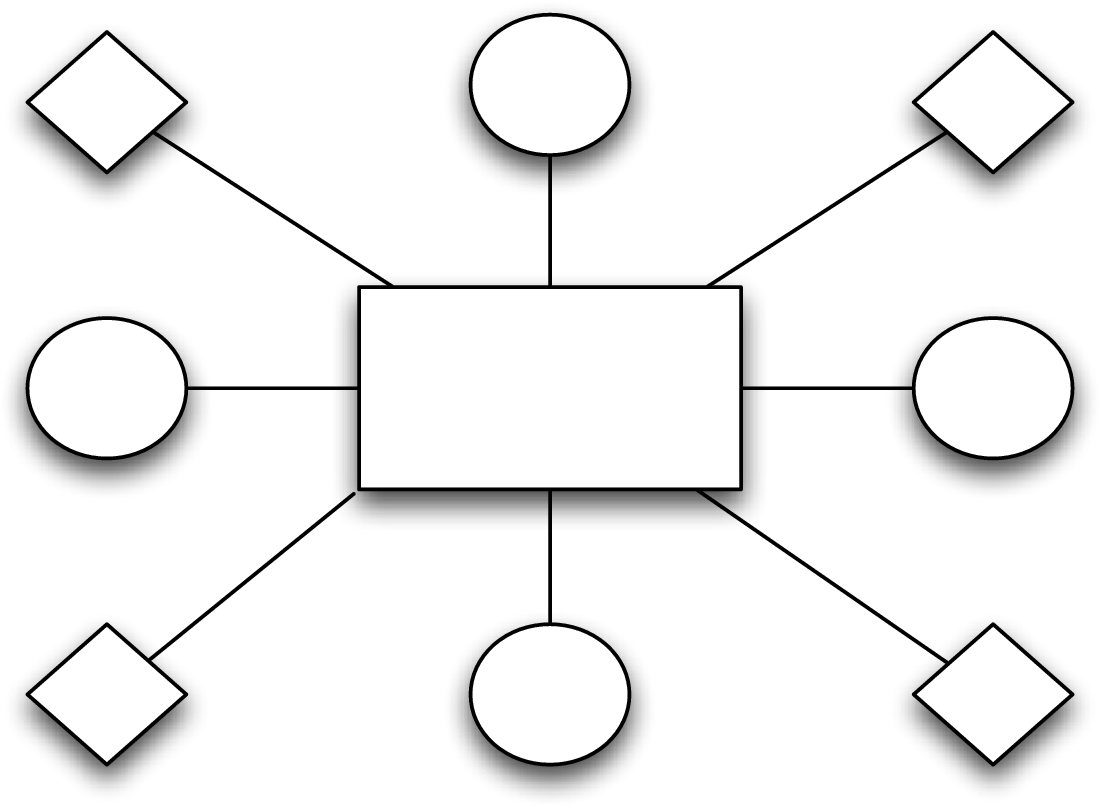}
  \caption{A schematic representation of networks analyzed in Section \ref{sec:basic_model} (left panel) and Section \ref{sec:basic_model_multiple} (right panel). The nodes are denoted with circles have a different dynamics from those indicated with squares. The dynamics of the common media is denoted with a rectangle.}
  \label{fig:quorum_1}
  \end{center}
\end{figure}

We now generalize the mathematical model analyzed in the previous Section, by allowing for $s \le N$ groups (or clusters) of nodes characterized by different dynamics (with possibly different dimensions) to communicate over the same common medium (see Figure \ref{fig:quorum_1}, right panel). We will prove a sufficient condition for the global exponential convergence of all nodes trajectories belonging to the same group towards each other. This regime is called concurrent synchronization~\cite{Pha_Slo_07}.

The mathematical model analyzed here is
\begin{equation}\label{eqn:general_model_non_homogeneous}
\begin{array}{*{20}l}
\dot x_i = f_{\gamma \left(i\right)}\left(x_i,z,t\right)\\
\dot z = g\left(z,\Psi\left(x_1,\ldots,x_N\right),t\right)\\ 
\end{array}
\end{equation}
where: i)  $\gamma$ is defined as in (\ref{eqn:gamma}); ii) $ x_i$ denotes the state variables of the network nodes (nodes belonging to different clusters may have different dimensions, say $n_{\gamma(i)}$) and $z$ denotes the state variables for the common medium ($z\in \R^d$); iii) $\Psi$, defined analogously to the previous Section, denotes the coupling function of the cluster $\gamma(i)$ with the common medium dynamics ($\Psi: \R^{n_{\gamma(1)}}\times \ldots\times \R^{n_{\gamma(N)}}\rightarrow \R^d$).

\begin{Theorem}\label{thm:quorum_convergence_non_homogeneous}
Concurrent synchronization is achieved in network (\ref{eqn:general_model_non_homogeneous}) if the functions $f_{\gamma\left(i\right)}\left(x,v(t),t\right)$ are all contracting for any $v(t) \in \R^d$.
\end{Theorem}
\proof Recall that (\ref{eqn:general_model_non_homogeneous}) is
composed by $N$ nodes having dynamics $f_1,\ldots,f_s$. Now, in analogy with the proof of Theorem \ref{thm:quorum_convergence}, consider the following virtual system:
\begin{equation} \label{eqn:virtual_quorum_non_homogeneous}
\begin{array}{*{20}l}
\dot y_1 = f_1\left(y_1,z,t\right)\\
\dot y_2 = f_2\left(y_2,z,t\right)\\
\vdots \\
\dot y_s = f_s\left(y_s,z,t\right)\\
\end{array}
\end{equation}
where $z(t)$ is seen as an exogenous input to the virtual system.
Let $\left\{X_{i}\right\}$ be the set of state variables belonging to
the $i$-th cluster composing the network, and denote with $X_{i,j}$ any
element of $\left\{X_i\right\}$. We have that
$(X_{1,j},\ldots,X_{s,j})$ are particular solutions of the virtual
system. Now, contraction of the virtual system implies that all of its
particular solutions converge towards each other, which in turn implies
that all the elements within the same cluster $\left\{X_i\right\}$
converge towards each other. Thus, contraction of the virtual system
(\ref{eqn:virtual_quorum_non_homogeneous}) implies concurrent
synchronization of the real system (\ref{eqn:general_model_non_homogeneous}).

To prove contraction of (\ref{eqn:virtual_quorum_non_homogeneous}), compute its Jacobian,
$$
J = \left[\begin{array}{*{20}c}
\frac{\partial f_1\left(y_1,z,t\right)}{\partial y_1} & 0 & 0  & \ldots & 0 \\
0 & \frac{\partial f_2\left(y_2,z,t\right)}{\partial y_2} & 0 & \ldots & 0  \\ 
\ldots & \ldots & \ldots & \ldots & \ldots\\
0 & 0 & 0 & 0 & \frac{\partial f_s\left(y_s,z,t\right)}{\partial y_s} \\
\end{array}\right]
$$
Now, by hypotheses, we have that all the functions $f_i(x,v(t),t)$ are contracting for any exogenous input. This in turn implies that the virtual system (\ref{eqn:virtual_quorum_non_homogeneous}) is contracting, since its Jacobian matrix is block diagonal with diagonal blocks being contracting.
 \endproof

\subsection{Systems communicating over different media}\label{sec:cluster}

\begin{figure}[thbp]
\begin{center}
  \includegraphics[width=8cm]{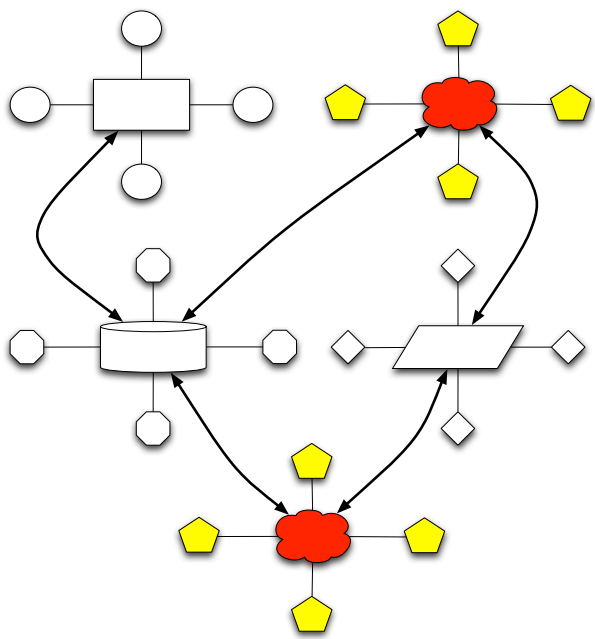}
  \caption{A schematic representation of the network analyzed Section \ref{sec:cluster}. The connections between media (and hence the connections of the autonomous level) are pointed out. Notice that only two nodes of the autonomous level are input equivalent since: i) their media have the same dynamics (in red); ii) both media are shared by the same number of nodes, sharing the same dynamics (in yellow).}
  \label{fig:quorum_2}
  \end{center}
\end{figure}

In the previous Section, we considered networks where some (possibly heterogeneous) nodes communicate over a common medium. We now consider a \emph{distributed} version of such topology, where each of the $s\le N$ groups composing the network have a private medium. Communication between the groups is then obtained by coupling only their media (see Figure \ref{fig:quorum_2}). The objective of this Section, is to provide a sufficient condition ensuring (concurrent) synchronization of such network topology. 

Note that the network topology considered here presents a layer structure. In analogy with the terminology used for describing the topology of the Internet and World-Wide-Web (see e.g. \cite{Boc_Lat_Mor_Cha_Hwa_06}, \cite{New_2003}), we term as \emph{medium} (or private) level the layer consisting of the nodes of the network and their corresponding (private) media; we then term as \emph{autonomous} level, the layer of the interconnections between the media. That is, the autonomous level is an \emph{abstraction} of the network, where its nodes' dynamics consists of the network nodes and their private medium. This in turn implies that in order for two nodes of the autonomous level to be identical they have to share: i) the same dynamics and number of nodes; ii) the same medium dynamics (see Figure \ref{fig:quorum_2}).

In what follows we will denote with $\mathcal{G}_p$ the set of homogeneous nodes communicating over the medium $z_p$. We will denote with $N_{p}$ the set of media which are linked to the medium $z_p$. Each medium communicates with its neighbors diffusively. The mathematical model is then:
\begin{equation}\label{eqn:general_model_groups}
\begin{array}{*{20}l}
\dot x_ i = f_p\left(x_i,z_p,t\right) & \ \ \ \ \ x_i \in \mathcal{G}_p\\
\dot z_p = g_p\left(z_p,\Psi\left(X_p\right),t\right)+ \sum_{j \in N_{p}} \left[\phi_p\left( z_j\right)-\phi_p\left( z_p\right) \right] & \ \ \ \ \ x_i \in \mathcal{G}_p
\end{array}
\end{equation}
where $p=1,\ldots,s$ and $X_p$ is the stack of all the vectors $x_i \in \mathcal{G}_p$. We assume that the dynamical equations for the media have all the same dimensions (e.g. $z_p\in R^d$), while the nodes belonging to different groups can have different dimensions (e.g. $x_i\in\R^p$, for any $i \in \mathcal{G}_p$). Here, the coupling functions between the media, $\phi_p:\R^d\rightarrow\R^d$, are assumed to be continuous and to have a diagonal Jacobian matrix with diagonal elements being nonnegative and bounded. All the matrices $\partial f_p/ \partial z$ are assumed to be bounded.

\begin{Theorem}\label{thm:quorum_convergence_groups}
Concurrent synchronization is attained in network (\ref{eqn:general_model_groups}) if: i) the nodes of its \emph{autonomous} level sharing the same dynamics are input equivalent; ii) $f_{p}\left(x_i,v(t),t\right)$,  $g_p\left(z_p,v(t),t\right)$ are all contracting functions for any $v(t) \in \R^d$; iii) $\frac{\partial f_p}{\partial z_p}$ are all uniformly bounded matrices.
\end{Theorem}
\proof
 Consider the following $2s$-dimensional virtual system, analogous to the one used for proving Theorem \ref{thm:quorum_convergence_non_homogeneous}:
\begin{equation}\label{eqn:virtual_general_model_groups}
\begin{array}{*{20}l}
\dot y_ {1,p} = f_p\left(y_{1,p},y_{2,p},t\right) \\
\dot y_{2,p} = g_p\left(y_{2,p},v_p(t),t\right)+ \sum_{k \in N_{p}} \left[\phi_p\left( y_{2,k}\right)-\phi_p\left( y_{2,p}\right) \right]
\end{array}
\end{equation}
where $p=1,\ldots,s$, and $v_p(t):=\Psi\left(X_p\right)$.
Notice that the above system is constructed in a similar way as (\ref{eqn:virtual_quorum_non_homogeneous}). In particular, solutions of (\ref{eqn:general_model_groups}) are particular solutions of the above virtual system (see the proof of Theorem \ref{thm:quorum_convergence_non_homogeneous}). That is, if cluster synchronization is attained for (\ref{eqn:virtual_general_model_groups}), then all the nodes sharing the same dynamics will converge towards each other. Now, Theorem \ref{thm:network_convergence_multiple_nodes} implies that cluster synchronization is attained for system (\ref{eqn:virtual_general_model_groups}) if: i) its nodes are contracting; ii) the coupling functions have a nonnegative bounded diagonal Jacobian; iii) nodes sharing the same dynamics are input equivalent. Since the last two conditions are satisfied by hypotheses, we have only to prove contraction of the (virtual) network nodes. In this view, differentiation of nodes dynamics in (\ref{eqn:virtual_general_model_groups}) yields the Jacobian matrix
$$
\left[\begin{array}{*{20}c}
\frac{\partial f_p\left(y_{1,p},y_{2,p},t\right)}{\partial y_{1,p}} & \frac{\partial f_p\left(y_{1,p},y_{2,p},t\right)}{\partial y_{2,p}}\\
0 & \frac{\partial g_p\left(y_{1,p},v_i(t),t\right)}{\partial y_{2,p}}
\end{array}\right]
$$
The above Jacobian has the structure of a hierarchy. Thus (see Section \ref{sec:math_tools}) the virtual system is contracting if:
\begin{enumerate}
\item $\frac{\partial f_p\left(y_{1,p},y_{2,p},t\right)}{\partial y_{1,p}}$ and $ \frac{\partial g_p\left(y_{2,p},v_i(t),t\right)}{\partial y_{2,p}}$ are both contracting
\item $\frac{\partial f_p\left(y_{1,p},y_{2,p},t\right)}{\partial y_{2,p}}$ is bounded
\end{enumerate}
The above two conditions are satisfied by hypotheses. Thus, the virtual network achieves cluster synchronization (Theorem \ref{thm:network_convergence_multiple_nodes}). This proves the Theorem.
 \endproof
 Note that Theorems \ref{thm:quorum_convergence} and
 \ref{thm:quorum_convergence_non_homogeneous} do not make any
 hypotheses on the medium dynamics $-$ synchronization (or concurrent
 synchronization) can be attained by the network nodes independently
 of the particular dynamics of the single medium, provided that the
 function $f$ (or the $f_i$'s) is contracting. By contrast, Theorem
 \ref{thm:quorum_convergence_groups} shows that the media dynamics
 becomes a key element for achieving concurrent synchronization in
 networks where different groups communicate over different media.

\section{Synchronization control}\label{sec:synchro_control}

In Section \ref{sec:main_res}, we derived several criteria ensuring
node synchronization for networks where
multiple nodes exchange their state variables using (multiple)
media. The above results also allow dimensionality reduction in the analysis
of the system's final behavior by treating each cluster as a single
element, similarly to \cite{Chu_Slo_Mil_07}, a point we will further illustrate
in Section~\ref{N}.

The objective of this Section is to provide a sufficient
condition guaranteeing some desired periodic behavior for the network
nodes. Specifically, we will guarantee a desired period for the steady
state oscillations. A related problem has been recently addressed in
\cite{Rus_diB_Son_09}, where entrainment of individual biological
systems to periodic inputs was analyzed.  We will now show
the following result,
\begin{Theorem}\label{thm:control_basic}
Consider the following network
\begin{equation}\label{eqn:general_mode_controll}
\begin{array}{*{20}l}
\dot x_ i = f\left(x_i,z,t\right) & \ \ \ \ \ i =1, \ldots, N\\
\dot z = g\left(z,\Psi\left(x_1,\ldots,x_N\right),t\right) + r\left(t\right)
\end{array}
\end{equation}
where $r\left(t\right)$ is a $T$-periodic signal.  All the nodes of the network synchronize onto a periodic orbit of period  $T$ if: i) $f\left(x_i,v(t),t\right)$ and $g\left(z,v(t),t\right)$ are contracting functions for any $v(t)\in \R^d$; ii) $\frac{\partial f}{\partial z}$ is bounded.
\end{Theorem}
\proof
Consider the following virtual system:
\begin{equation}\label{eqn:virtual_quorum_controlled}
\begin{array}{*{20}l}
\dot y_1 = f\left(y_1,y_2,t\right)\\
\dot y_2 = g\left(y_2,v(t),t\right) + r\left(t\right)\\
\end{array}
\end{equation}
where $v(t):=\Psi\left(x_1,\ldots,x_N\right)$. We will prove the Theorem by showing that such a system is contracting. Indeed, in this case, the trajectories of (\ref{eqn:virtual_quorum_controlled}) will globally exponentially converge to a unique $T$-periodic solution, implying that also $x_i$ will exhibit a $T$-periodic steady state behavior. Differentiation of the virtual system yields:
$$
\left[\begin{array}{*{20}c}
\frac{\partial f\left(y_1,y_2,t\right)}{\partial y_1} & \frac{\partial f\left(y_1,y_2,t\right)}{\partial y_2}\\
0 & \frac{\partial g\left(y_2,v(t),t\right)}{\partial y_2}
\end{array}\right]
$$
The above Jacobian has the structure of a hierarchy. Thus (see Section \ref{sec:math_tools}) the virtual system is contracting if:
\begin{enumerate}
\item $\frac{\partial f\left(y_1,y_2,t\right)}{\partial y_1}$ and $ \frac{\partial g\left(y_2,v(t),t\right)}{\partial y_2}$ are both contracting
\item $\frac{\partial f\left(y_1,y_2,t\right)}{\partial y_2}$ is bounded
\end{enumerate}
The first condition is satisfied since, by hypotheses, the functions
$f\left(x,v(t),t\right)$ and $g\left(z,v(t),t\right)$ are
contracting for any $v\in \R^d$. The second condition is also
satisfied since we assumed $\partial f/\partial z$ to be bounded. The Theorem is then proved.
\endproof

The results can be extended to the more general case of networks of
non homogeneous nodes communicating over non homogeneous media.

\begin{Theorem}\label{thm:control_concurrent}
Consider the following network
\begin{equation}\label{eqn:general_model_groups_control}
\begin{array}{*{20}l}
\dot x_ i = f_p\left(x_i,z_p,t\right) & \ \ \ \ \ x_i \in \mathcal{G}_p\\
\dot z_p = g_p\left(z_p,\Psi\left(X_p\right),t\right) + \sum_{k \in N_{p}} \left[\phi\left( z_k\right)-\phi\left( z_p\right) \right]+r\left(t\right) & \ \ \ \ \ x_j \in \mathcal{G}_p
\end{array}
\end{equation}
where $X_p$ is the stack of all the $x_i \in \mathcal{G}_p$ and $r\left(t\right)$ is a $T$-periodic signal. Concurrent synchronization is attained, with a steady state periodic behavior of period $T$ if:
\begin{enumerate}
\item the nodes of the autonomous level sharing the same dynamics are input equivalent;
\item the coupling functions $\phi$ have bounded diagonal Jacobian with nonnegative diagonal elements;
\item $f_{p}\left(x_i,v(t),t\right)$ and $g_p\left(z_p,v(t),t\right)$ are contracting functions for any $v(t) \in \R^d$;
\item $\frac{\partial f_p}{\partial z_p}$ are all uniformly bounded matrices.
\end{enumerate}
\end{Theorem}
\proof The proof is formally the same as that of Theorem \ref{thm:quorum_convergence_groups} and Theorem \ref{thm:control_basic},
and it is omitted here for the sake of brevity.
\endproof

\subsubsection*{A simple example}

Consider a simple biochemical reaction, consisting of
a set of $N>1$ enzymes sharing the same substrate. We
denote with $X_1,\ldots, X_N$ the
concentration of the reaction products. We also assume that the
dynamics of $S$ is affected by some $T$-periodic input, $r(t)$ (the behavior of networks where the medium dynamics is affected by an exogenous input will be analyzed in Section \ref{sec:synchro_control}). We
assume that the total concentration of $X_{i}$, i.e. $X_{i,T}$, is
much less than the initial substrate concentration, $S_0$. In these
hypotheses, a suitable mathematical model for the system is given by
(see e.g. \cite{Sza_Ste_Per_06}):
\begin{equation}\label{eqn:simple_example_model}
\begin{array}{*{20}l}
\dot X_i = - a X_i + \frac{K_{1} S}{K_{2} +S} & i=1,\ldots, N\\
\dot S = - \sum_{i=1}^N \frac{K_{1} S}{K_{2} +S} + r(t) \\
\end{array}
\end{equation}
with $K_1$ and $K_2$ be positive parameters.
Thus, a suitable virtual system for the network is
\begin{equation}\label{eqn:simple_example_virtual}
\begin{array}{*{20}l}
\dot y_1 = -a y_1 + \frac{K_{1} y_2}{K_{2} +y_2}\\
\dot y_2 = - \sum_{i=1}^N \frac{K_{1} y_2}{K_{2} +y_2} + r(t) \\
\end{array}
\end{equation}
Differentiation of the above system yields the Jacobian matrix
\begin{equation}\label{eqn:simple_example_Jacobian}
\left[\begin{array}{*{20}c}
-a & \frac{K_2}{(K_2+y_1)^2}\\
0 & -N\frac{K_2}{(K_2+y_1)^2}\\
\end{array}\right]
\end{equation}
It is straightforward to check that the above matrix represents a contracting hierarchy. Thus, all the trajectories of the virtual system globally exponentially converge towards a unique $T$-periodic solution. This, in turn, implies that $X_i$, $i=1,\ldots, N$, globally exponentially converge towards each other and towards the same periodic solution.

\begin{figure}[thbp]
\begin{center}
\centering \psfrag{x}[c]{{time (minutes)}}
\centering \psfrag{y3}[c]{{r(t)}}
\centering \psfrag{y1}[c]{{$X_i(t)$}}
\centering \psfrag{y2}[c]{{$S(t)$}}
    \includegraphics[width=10cm]{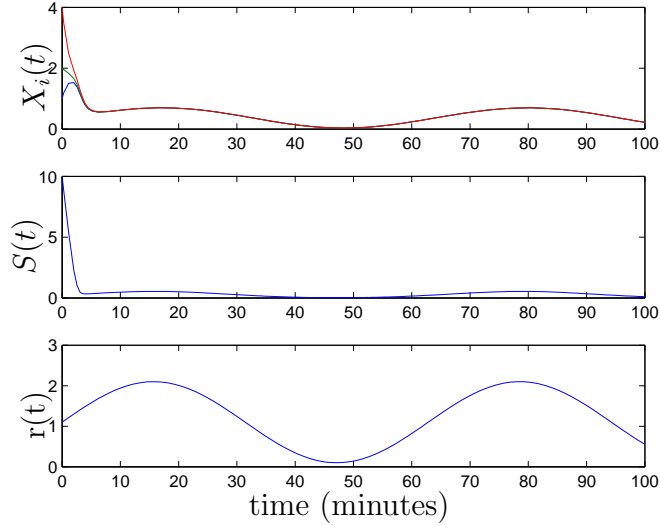}
  \caption{Simulation of (\ref{eqn:simple_example_model}), with $N=3$ and $r(t) = 1.1 + \sin(0.1*t)$. System parameters are set as follows: $a=1$, $K_2=1$, $K_1=2$.}
  \label{fig:quorum_enzyme}
  \end{center}
\end{figure}

Figure \ref{fig:quorum_enzyme} illustrates the behavior for $N=3$. Notice that, as expected from the above theoretical analysis, $X_1$, $X_2$ and $X_3$ synchronize onto a periodic orbit of the same period as $r(t)$. 

\section{Emergent properties as $N$ increases}\label{N}

In this Section, we analyze how the convergence properties of  
a given quorum sensing network vary as the number $N$ of nodes increases.
We show that for typical quorum sensing networks, as 
$N$ becomes sufficiently large, synchronization always occurs.
One particular modeling context where these results have important
implications is that of cell proliferation in biological systems.

\subsection{A lower bound on $N$ ensuring synchronization}

It is well known~\cite{Wan_Slo_05} that for all-to-all diffusively
coupled networks of the form
\begin{equation}\label{eqn:all_to_all_no_dyn}
\dot x_i = f(x_i,t) + \sum_{i=1}^Nk(x_j-x_i)
\end{equation}
the minimum coupling gain $k$ required for synchronization is inversely proportional to the number of nodes composing the network. That is,
$$
k_{\min} \propto \frac{1}{N}
$$ 
We now show that a similar bound holds for nodes coupled by means of quorum sensing of the form
\begin{equation}\label{eqn:general_model_part_simplified}
\begin{array}{*{20}l}
\dot x_ i = f\left(x_i,t\right) + kN (z -  x_i) & \ \ \ \ \ i =1, \ldots, N\\
\dot z = g\left(z,\Psi(x_1,\ldots,x_N),t\right)
\end{array}
\end{equation}
To simplify notations, the above model assumes that $z$ and all $x_i$ have the same dimensions. Also note that in (\ref{eqn:general_model_part_simplified}) the dependence of the coupling gain on the number of nodes, $N$, is given explicitly.
\begin{Theorem}\label{thm:threshold_condition_sync}
Assume that the Jacobian $\left(\frac{\partial f}{\partial x}\right)$ is upper-bounded
by $\alpha$ for some matrix measure $\mu$, i.e.,
$$
\exists \alpha \in \R,\ \forall x,\ \forall t \ge 0, \ \ \ \ \ \mu\left(\frac{\partial f}{\partial x}\right) \le \alpha
$$
Then, network (\ref{eqn:general_model_part_simplified}) synchronizes if
$$
k>\frac{\alpha}{N}
$$
That is, $k_{\min}\propto 1/N \ $.
\end{Theorem}
\proof
Consider the virtual system
\begin{equation}\label{eqn:proof_measure}
\dot y = f(y,t) + kN(z - y)
\end{equation}
Synchronization is attained if the virtual system is contracting.
Now, computing the matrix measure of the Jacobian of
(\ref{eqn:proof_measure}) yields
$$
 \forall x,\ \forall t \ge 0, \ \ \ \ \ \ \ \ \mu\left(\frac{\partial f}{\partial y} - kN I\right)
\ \le \ \mu\left(\frac{\partial f}{\partial y}\right) + kN \mu\left(-  I\right)
\ \le \ \alpha - kN
$$
Thus, the virtual system is contracting if $k>\frac{\alpha}{N}$.
\endproof

\subsection{Dependence on initial conditions}

We now consider the basic quorum sensing model
(\ref{eqn:general_model}). We derive simple conditions for the
final behavior of the network to become independent of initial
conditions (in the nodes and the medium) as $N$ becomes large.

\begin{Theorem}\label{thm:quorum_convergence_global}
Assume that for (\ref{eqn:general_model}) the following conditions hold: 
\begin{itemize}
\item $\mu\left(\frac{\partial f}{\partial x}\right)\rightarrow - \infty$ as $N\rightarrow +\infty$ 
\item $g\left(z,v_2(t),t\right)$ is contracting (for any $v_2(t)$ in $\R^d$)
 \item $\norm{\frac{\partial f}{\partial z}}$ and $\norm{\frac{\partial g}{\partial v_2}}$ are bounded for any $x$, $z$, $v_2$ (where  $\norm{\cdot}$ is the operator norm)

\end{itemize}
Then, there exists some $N^\ast$ such that for any $N\ge N^\ast$ all
  trajectories of (\ref{eqn:general_model}) globally exponentially
  converge towards a unique synchronized solution, independent of
  initial conditions.
\end{Theorem}

\proof We know that contraction of $f\left(x,v_1(t),t\right)$ for any
$v_1(t)$ (which the first condition implies for $N$ large enough)
ensures network synchronization. That is, there exists a unique
trajectory, $x_s(t)$, such that, as $t\rightarrow + \infty$,
$$
\abs{x_i-x_s} \rightarrow 0, \quad \forall i
$$
Therefore, the final behavior is described by the following lower-dimensional system:
\begin{equation}\label{eqn:general_model_reduced}
\begin{array}{*{20}l}
\dot x_ s = f\left(x_s,z,t\right)\\
\dot z = g\left(z, \Psi\left(x_s\right),t\right)
\end{array}
\end{equation}
If in turn this reduced-order system (\ref{eqn:general_model_reduced})
is contracting, then its trajectories globally exponentially
converge towards a unique solution, say $x_s^\ast(t)$, regardless of
initial conditions. This will prove the Theorem (similar strategies
are extensively discussed in \cite{Chu_Slo_Mil_07}).

To show that (\ref{eqn:general_model_reduced}) is indeed contracting,
compute its Jacobian matrix,
\begin{equation}\nonumber
\left[\begin{array}{*{20}c}
\frac{\partial f}{\partial x_s} & \frac{\partial f}{\partial z}\\
\frac{\partial g}{\partial x_s} & \frac{\partial g}{\partial z}\\
\end{array}\right]
\end{equation}
Lemma \ref{lem:partition} in the Appendix shows that the above matrix is contracting if there exists some strictly positive constants $\theta_1$, $\theta_2$ such that
\begin{equation}\label{eqn:cond_global_sync}
\mu\left(\frac{\partial f}{\partial x_s}\right) + \frac{\theta_2}{\theta_1}\norm{\frac{\partial g}{\partial x_s}}\ \ \ \ \ \rm{and}\ \ \ \ \ 
\mu\left(\frac{\partial g}{\partial z}\right) + \frac{\theta_1}{\theta_2}\norm{\frac{\partial f}{\partial z}}
\end{equation}
are both uniformly negative definite.

Now, $\mu\left(\frac{\partial f}{\partial x_s}\right)$ and
$\mu\left(\frac{\partial g}{\partial z}\right)$ are both uniformly
negative by hypotheses. Furthermore, $\mu\left(\frac{\partial
    f}{\partial x_s}\right)$ tends to $-\infty$ as $N$ increases: since
$\norm{\frac{\partial f}{\partial z}}$ and $\norm{\frac{\partial
    g}{\partial x_s}}$ are bounded, this implies that there exists
some $N^\ast$ such that for any $N\ge N^\ast$ the two conditions in
(\ref{eqn:cond_global_sync}) are satisfied.
 \endproof
 
 Also, assume that actually the dynamics $f$ and $g$ do not depend
 explicitly on time. Then, under the conditions of the above Theorem,
 the reduced system is both contracting and autonomous, and so it tends
 towards a unique equilibrium point~\cite{Loh_Slo_98}. Thus, the
 original system converges to a unique equilibrium, where all $x_i$'s
 are equal.

 In addition, note that when the synchronization rate and the
 contraction rate of the reduced system both increase with $N$, this
 also increases robustness~\cite{Pha_Slo_07} to variability and
 disturbances.

\subsection{How synchronization protects from noise}

In this section, we discuss briefly how the synchronization mechanism
provided by dynamical quorum sensing protects from noise and
variability in a fashion similar to the static mechanism studied in
\cite{Tab_Slo_Pha_09}. We show that the results
of~\cite{Tab_Slo_Pha_09}, to which the reader is referred for details
about stochastic tools, extend straightforwardly to the case where the
dynamics of the quorum variables cannot be neglected or indeed may
play a central part, as studied in this paper.

Assume that the dynamics of each network element $x_i$ in
(\ref{eqn:general_model_part_simplified}) is subject to noise, and
consider, similarly to~\cite{Tab_Slo_Pha_09}, the corresponding system
of individual elements in Ito form
\begin{equation} \label{eq:main}
  dx_i= \left(f(x_i,t) + kN(z - x_i)\right)dt + \sigma dW_i\ \ \ \ \ \ \ \ \ i=1\dots N
\end{equation}
where the all-to-all coupling in ~\cite{Tab_Slo_Pha_09} has been
replaced by a more general quorum sensing mechanism.  The subsystems
are driven by independent noise processes, and for simplicity the
noise intensity $\sigma$ in the equations above is assumed to be
constant. We make no assumptions about noise acting directly on the
dynamics of the environment/quorum vector $z$.

Proceeding exactly as in \cite{Tab_Slo_Pha_09} yields similar results
on the effect of noise. In particular, let $x^\bullet$ be the center
of mass of the $x_i$, that is
\begin{equation}\nonumber
  x^\bullet= \frac{1}{N} \sum_i x_i \
\end{equation}

Notice that when all the nodes are synchronized onto some common solution, say 
$x_s(t)$, then, by definition, $x^\bullet =x_s(t)$.

Adding up the dynamics in (\ref{eq:main}) gives
\begin{equation} \label{eq:xmass-dynamic}
dx^\bullet = \frac{1}{N}\left(\sum_{i} f(x_i,t)\right) dt + kN(z-x^\bullet)dt
+ \frac{1}{N}\sum_i \sigma dW_i
\end{equation}
Let
\begin{equation}\nonumber
  \epsilon= f(x^\bullet,t) - \frac{1}{N}\left(\sum_{i=1}^N f(x_i,t)\right)
\end{equation}
Note that $\epsilon=0$ when all the nodes are synchronized.
        
By analogy with (\ref{eq:main}), equation (\ref{eq:xmass-dynamic}) can then be written
\begin{equation} \label{average}
  dx^\bullet = \left(f(x^\bullet,t) + kN(z- x^\bullet)  + \epsilon \right) dt +
  \frac{1}{N}\sum_i \sigma dW_i
\end{equation}
Using the Taylor formula with integral remainder exactly as in \cite{Tab_Slo_Pha_09} yields
a bound on the distortion term $\epsilon$, as a function of the nonlinearity, the coupling\
 gain $k$, and the number of cells $N$,
\begin{equation}\nonumber
\E(\|\epsilon\|) \ \le \ \lambda_{\rm max}\left(\frac{\partial^2 f}{\partial x^2}\right)\ \ \rho(kN)
\end{equation}
\noindent where $\lambda_{\rm max}(\frac{\partial^2 f}{\partial x^2})$
is a uniform upper bound on the spectral radius of the Hessian
$\frac{\partial^2 f}{\partial x^2}$, and $\rho(kN) \rightarrow 0$ as
$kN \rightarrow +\infty$.  In particular, in (\ref{average}),
both the distortion term $\epsilon$ and the average noise term $\
\frac{1}{N}\sum_i \sigma dW_i\ $ tend to zero as $N \rightarrow + \infty$.

Note that an additional source of noise may be provided by the
environment on the quorum variables themselves.  We made no
assumptions above about such noise, which acts directly on the
dynamics of the environment/quorum vector $z$. How it specifically
affects the common quantity $z$ in (\ref{eq:main}) could be further
studied.

Similar results hold for the effects of bounded disturbances and dynamic variations.

\section{Examples}\label{sec:examples}

\subsection{Controlling synchronization of genetic relaxation oscillators}\label{sec:synchro_example}

\begin{figure}[thbp]
\begin{center}
  \includegraphics[width=8cm]{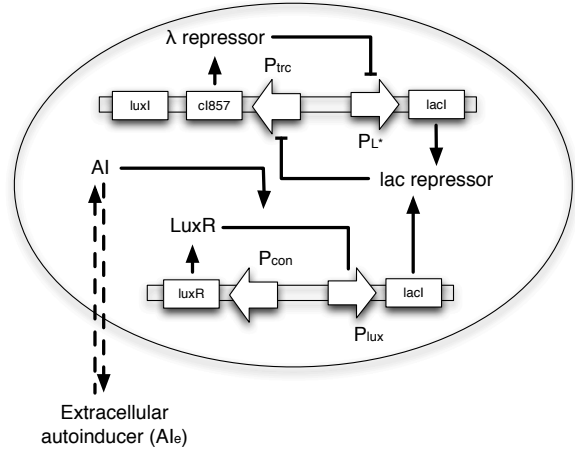}
  \caption{A schematic representation of the genetic circuit: detailed circuit.}
  \label{fig:genetic_oscillator}
  \end{center}
\end{figure}

We now consider the problem of synchronizing a population of genetic oscillators. Specifically, we consider the genetic circuit analyzed in \cite{Kuz_Kae_Kop_04} (a variant of \cite{Kob_Kae_Ara_Chu_Gar_Can_Col_04}), and schematically represented in Figure \ref{fig:genetic_oscillator}. Such a circuit is composed of two engineered gene networks that have been experimentally implemented in \emph{E. coli}; namely: the toggle switch \cite{Gar_Can_Col_00} and an intercell communication system \cite{You_Cox_Wei_Arn_04}. The toggle switch is composed of two transcription factors: the \emph{lac} repressor, encoded by gene \emph{lacI}, and the temperature-sensitive variant of the $\lambda cI$ repressor, encoded by the gene \emph{cI857}. The expressions of \emph{cI8547} and \emph{lacI} are controlled by the promoters $P_{trc}$ and $P_{L^\ast}$ respectively (for further details see \cite{Kuz_Kae_Kop_04}). The intercell communication system makes use of components of the quorum-sensing system from \emph{Vibro fischeri} (see e.g. \cite{Ng_Bas_09} and references therein). Such a mechanism allows cells to sense population density through the transcription factor LuxR, which is an activator of the genes expressed by the $P_{lux}$ promoter, when a small molecule $AI$ binds to it. This small molecule, synthesized by the protein LuxI, is termed as autoinducer and it can diffuse across the cell membrane.  

\begin{figure}[thbp]
\begin{center}
  \includegraphics[width=8cm]{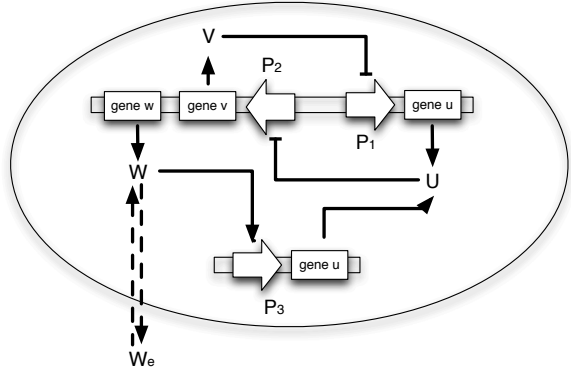}
  \caption{Simplified circuit using for deriving the mathematical model (\ref{eqn:math_model_relax_osc}). Both the promoters and transcription factors are renamed.}
  \label{fig:genetic_oscillator_simplified}
  \end{center}
\end{figure}

In \cite{Kuz_Kae_Kop_04}, the following dimensionless simplified model is analyzed (see Figure \ref{fig:genetic_oscillator_simplified}):
\begin{subequations}\label{eqn:math_model_relax_osc}
\begin{equation}\label{eqn:math_model_relax_osc_1}
\dot u_i = \frac{\alpha_1}{1+v_i^\beta} + \frac{\alpha_3 w_i^\eta}{1+w_i^\eta} -d_1 u_i
\end{equation}
\begin{equation}\label{eqn:math_model_relax_osc_2}
\dot v_i = \frac{\alpha_2}{1+u_i^\gamma} - d_2 v_i
\end{equation}
\begin{equation}\label{eqn:math_model_relax_osc_3}
\dot w_i = \varepsilon \left(\frac{\alpha_4}{1+u_i^\gamma} - d_3w_i\right)+2 d\left(w_e-w_i\right)
\end{equation}
\begin{equation}\label{eqn:math_model_relax_osc_4}
\dot w_e = \frac{D_e}{N}\sum_{i=1}^N\left(w_i-w_e\right)-d_e w_e
\end{equation}
\end{subequations}
where $u_i$, $v_i$ and $w_i$ denotes the (dimensionless) concentrations of the \emph{lac} repressor, $\lambda$ repressor and LuxR-AI activator respectively. The state variable $w_e$ denotes instead the (dimensionless) concentration of the extracellular autoinducer.

In \cite{Kuz_Kae_Kop_04}, a bifurcation analysis is performed for the
above model, showing that synchronization can be attained for some
range of the biochemical parameters of the circuit. However, as the
objective of that paper was to analyze the onset of synchronization,
the problem of guaranteeing a desired oscillatory behavior was not
addressed. In what follows, using the results derived in the previous
sections, we address the open problem of guaranteeing a desired period
for the steady state oscillatory behavior of network
(\ref{eqn:math_model_relax_osc}).

\begin{figure}[thbp]
\begin{center}
  \includegraphics[width=8cm]{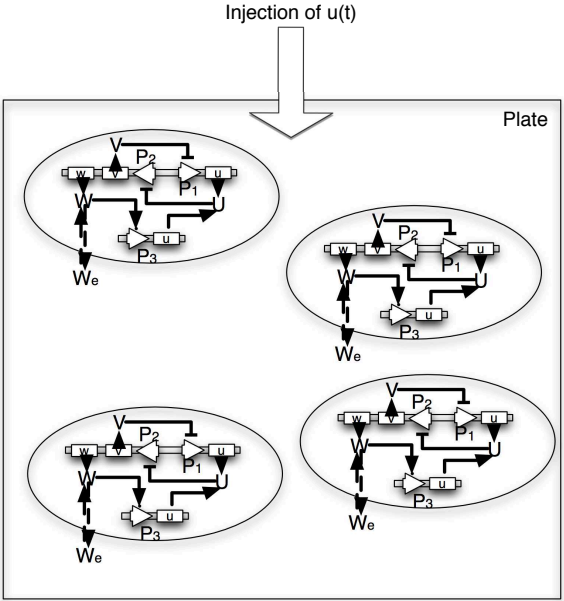}
  \caption{Network control setup.}
  \label{fig:genetic_oscillator_simplified_network}
  \end{center}
\end{figure}

The control mechanism that we use here is an exogenous signal acting on the extracellular autoinducer concentration, see also \cite{Rus_diB_Son_09}. That is, the idea is to modify (\ref{eqn:math_model_relax_osc_4}) as follows
\begin{equation}\label{eqn:math_model_relax_osc_4_mod}
\dot w_e = \frac{D_e}{N}\sum_{i=1}^N\left(w_i-w_e\right)-d_e w_e + r\left(t\right)
\end{equation}
where $r\left(t\right)$ is some $T$-periodic signal.
The set up that we have in mind here is illustrated in Figure \ref{fig:genetic_oscillator_simplified_network}, where multiple copies of the genetic circuit of interest share the same surrounding solution, on which $r\left(t\right)$ acts. From the technological viewpoint,  $r\left(t\right)$ can be implemented by controlling the temperature of the surrounding solution, and/or using e.g. the recently developed microfluidics technology (see e.g. \cite{Bee_Men_Wal_02} and references therein).

In what follows, we will use Theorem \ref{thm:quorum_convergence} to
find a set of biochemical parameters that ensure synchronization of
(\ref{eqn:math_model_relax_osc_1})-(\ref{eqn:math_model_relax_osc_4}). This,
using the results of Section~\ref{sec:synchro_control}, immediately implies that the
forced network
(\ref{eqn:math_model_relax_osc_1})-(\ref{eqn:math_model_relax_osc_3}),
(\ref{eqn:math_model_relax_osc_4_mod}) globally exponentially
converges towards a $T$-periodic steady state behavior.

System (\ref{eqn:math_model_relax_osc}) has the same structure as (\ref{eqn:general_model_part}), with $x_i =\left[u_i, v_i,w_i\right]^T$, $z = w_e$, and:
$$
\begin{array}{*{20}l}
f\left(x_i,t\right) = \left[\begin{array}{*{20}c}
\frac{\alpha_1}{1+v_i^\beta} + \frac{\alpha_3 w_i^\eta}{1+w_i^\eta} - d_1u_i\\
 \frac{\alpha_2}{1+u_i^\gamma} - d_2v_i\\
 \varepsilon \left(\frac{\alpha_4}{1+u_i^\gamma} - d_3w_i\right)\\
\end{array}\right] & k_z\left( z \right) - k_x\left( x_i\right) = \left[\begin{array}{*{20}c}
0\\
0\\
2 d\left(w_e-w_i\right)\\
\end{array}\right]\\
g\left(z,t\right) = -d_e w_e & \sum_{i=1}^N\left[u_x\left( x_i\right) - u_z \left( z \right)\right] = \frac{D_e}{N}\sum_{i=1}^N\left(w_i-w_e\right)
\end{array}
$$
We know from Theorem  \ref{thm:quorum_convergence} that all nodes trajectories converge towards each other if:
\begin{enumerate}
\item \ $f\left(x_i,t\right) - k_x\left( x_i\right)$ is contracting;
\item \ $g\left(z,t\right)- N u_z\left( z\right)$ is contracting.
\end{enumerate}
That is, contraction is ensured if there exist some matrix measures, $\mu_\ast$ and $\mu_{\ast\ast}$, such that
$$
\begin{array}{*{20}c}
\mu_\ast\left(\left(x_i,t\right) - k_x\left( x_i\right)\right) \ \ \ {\rm and}\ \ & \mu_{\ast\ast}\left(g\left(z,t\right)- N u_z\left( z\right)\right)
\end{array}
$$
are uniformly negative definite.
We use the above two conditions in order to obtain a set of
biochemical parameters ensuring node convergence. A possible choice
for the above matrix measures is $\mu_{\ast}=\mu_{\ast\ast} = \mu_1$
(see \cite{Rus_diB_Slo_09,Rus_diB_Son_09}). Clearly, other choices for
the matrix measures $\mu_{\ast}$ and $\mu_{\ast\ast}$ can be made,
leading to different algebraic conditions, and thus to (eventually)
a different choice of biochemical parameters.

We assume that $\beta =\eta=\gamma =2$, and show how to find a set of biochemical parameters satisfying the above two conditions.

{\bf Condition 1.} Differentiation of $\frac{\partial f}{\partial x_i} - \frac{\partial k}{\partial x_i}$ yields the Jacobian matrix (where the subscripts have been omitted)
\begin{equation}\label{eqn:Jacobian_relax_osc}
J_i := \left[\begin{array}{*{20}c}
-d_1 & \frac{-2\alpha_1v}{\left(1+v^2\right)^2} & \frac{2\alpha_3w}{\left(1+w^2\right)^2}\\
\frac{-2\alpha_2u}{\left(1+u^2\right)^2} & -d_2 & 0\\
\frac{-2\varepsilon\alpha_4u}{\left(1+u^2\right)^2} & 0 & -\varepsilon d_3-2d\\
\end{array}\right]
\end{equation}
Now, by definition of $\mu_1$, we have:
$$
\mu_1\left(J_i\right) = \max\left\{-d_1 +\frac{2\alpha_2u}{\left(1+u^2\right)^2}+ \frac{2\varepsilon\alpha_4u}{\left(1+u^2\right)^2}, -d_2 + \frac{2\alpha_1v}{\left(1+v^2\right)^2}, -\varepsilon d_3 -2d + \frac{2\alpha_3w}{\left(1+w^2\right)^2} \right\}
$$ 
Thus, $J_i$ is contracting if $\mu_1\left(J_i\right)$ is uniformly negative definite. That is,
\begin{equation}\label{eqn:condition_synchro}
\begin{array}{*{20}l}
-d_1 +\frac{2\alpha_2u}{\left(1+u^2\right)^2}+ \frac{2\varepsilon\alpha_4u}{\left(1+u^2\right)^2}\\
 -d_2 + \frac{2\alpha_1v}{\left(1+v^2\right)^2}\\
-\varepsilon d_3 -2d + \frac{2\alpha_3w}{\left(1+w^2\right)^2}\\
\end{array}
\end{equation}
are all uniformly negative. Notice now that the maximum of the function $a\left(v\right) = \frac{\bar a v}{\left(1+v^2\right)^2}$ is  $\hat a = \frac{3\sqrt{3} \bar a}{16}$. Thus, the set of inequalities (\ref{eqn:condition_synchro}) is fulfilled if:
\begin{equation}\label{eqn:condition_synchro_max}
\begin{array}{*{20}l}
-d_1 + \frac{6\alpha_2\sqrt{3}}{16} + \frac{6\varepsilon \alpha_4\sqrt{3}}{16}\\
 -d_2 + \frac{6\alpha_1\sqrt{3}}{16}\\
-\varepsilon d_3 -2d +\frac{6\alpha_3\sqrt{3}}{16}\\
\end{array}
\end{equation}
are all uniformly negative.

{\bf Condition 2} In this case it is easy to check that the matrix $J_e := \frac{\partial g}{\partial z}- N \frac{\partial u}{\partial z}$ is contracting for any choice of the (positive) biochemical parameters $D_e$, $d_e$.

\begin{figure}[thbp]
\begin{center}
\centering \psfrag{x1}[c]{{time (minutes)}}
\centering \psfrag{x2}[c]{{time (minutes)}}
\centering \psfrag{y1}[c]{{$w_i$ (arbitrary units)}}
\centering \psfrag{y2}[c]{{$r\left(t\right)$}}
  \includegraphics[width=8cm]{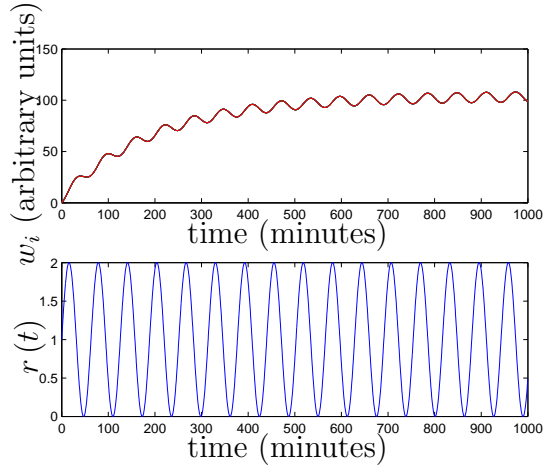}
  \caption{Behavior of (\ref{eqn:math_model_relax_osc_1})-(\ref{eqn:math_model_relax_osc_3}), (\ref{eqn:math_model_relax_osc_4_mod}), when forced by $r\left(t\right) = 1+\sin\left(0.1t\right)$. Notice that the nodes have initial different conditions, and that they all converge onto a common asymptotic having the same period as $r\left(t\right)$.}
  \label{fig: sim_synchro_genetic}
  \end{center}
\end{figure}

Thus, we can conclude that any choice of biochemical parameters fulfilling (\ref{eqn:condition_synchro_max}) ensures synchronization of the network onto a periodic orbit of period $T$.  In \cite{Kuz_Kae_Kop_04}, it was shown that a set of parameters for which synchronization is attained is: $\alpha_1=3$, $\alpha_2=4.5$, $\alpha_3=1$, $\alpha_4=4$, $\varepsilon = 0.01$, $d=2$, $d_1=d_2=d_3=1$. We now use the guidelines provided by (\ref{eqn:condition_synchro_max}) to make a minimal change of the parameters values ensuring network synchronization with steady state oscillations of period $T$. Specifically, such conditions can be satisfied by setting $d_1=6$, $d_2=2$.  Figure \ref{fig: sim_synchro_genetic} shows the behavior of the network for such a choice of the parameters. 

\subsubsection{Biological oscillators communicating over different media}

\begin{figure}[thbp]
\begin{center}
  \includegraphics[width=10cm]{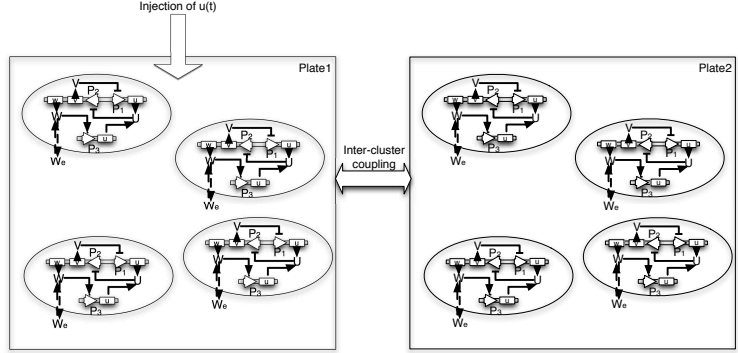}
  \caption{Two clusters of genetic circuits communicating over two different media}
  \label{fig:genetic_oscillator_simplified_network_cluster}
  \end{center}
\end{figure}

In the above Section, we assumed that all the genetic circuits shared the same surrounding solution. We now analyze the case where two different clusters of genetic circuits are surrounded by two different media. The communication between clusters is then left to some (eventually artificial) communication strategy between the two media (see Figure \ref{fig:genetic_oscillator_simplified_network_cluster}).

Notice that only one of the two media is forced by the exogenous $T$-periodic signal $r\left(t\right)$ (thus the dynamics of the two clusters are not the same), while the two media communicate with each other in a diffusive way. The mathematical model that we analyze here is then:
\begin{equation}\label{eqn:math_model_relax_osc_cluster}
\begin{array}{*{20}l}
\dot u_{i1} = \frac{\alpha_1}{1+v_{i1}^\beta} + \frac{\alpha_3 w_{i1}^\eta}{1+w_{i1}^\eta} -d_1 u_{i1}\\
\dot v_{i1} = \frac{\alpha_2}{1+u_{i1}^\gamma} - d_2 v_{i1}\\
\dot w_{i1} = \varepsilon \left(\frac{\alpha_4}{1+u_{i1}^\gamma} - d_3w_{i1}\right)+2 d\left(w_{e1}-w_{i1}\right)\\
\dot w_{e1} = \frac{D_{e}}{N}\sum_{i=1}^N\left(w_{i1}-w_{e1}\right)-d_{e} w_{e1}+ r\left(t\right) + \phi\left(w_{e2}\right) - \phi \left(w_{e1}\right)\\
\dot u_{i2} = \frac{\alpha_1}{1+v_{i2}^\beta} + \frac{\alpha_3 w_{i2}^\eta}{1+w_{i2}^\eta} -d_1 u_{i2}\\
\dot v_{i2} = \frac{\alpha_2}{1+u_{i2}^\gamma} - d_2 v_{i2}\\
\dot w_{i2} = \varepsilon \left(\frac{\alpha_4}{1+u_{i2}^\gamma} - d_3w_{i2}\right)+2 d\left(w_{e2}-w_{i2}\right)\\
\dot w_{e2} = \frac{D_{e}}{N}\sum_{i=1}^N\left(w_{i2}-w_{e2}\right)-d_{e} w_{e2}+  \phi\left(w_{e1}\right) - \phi \left(w_{e2}\right)\\
\end{array}
\end{equation}
where $x_{i1} =\left[u_{i1}, v_{i1}, w_{i1}\right]^T$ and $x_{i2} =\left[u_{i2}, v_{i2}, w_{i2}\right]^T$ denote the set of state variables of the $i$-th oscillator of the first and second cluster respectively. Analogously, $w_{e1}$ and $w_{e2}$ denote the extracellular autoinducer concentration surrounding the first and second cluster of genetic circuits. In the above model we assume that the biochemical parameters of the two genetic circuits and media are the same.

To ensure concurrent synchronization,  we  tune the biochemical parameters of the two clusters of oscillators and design the coupling function between the media ($\phi\left(\cdot\right)$) by using the guidelines provided by Theorem \ref{thm:quorum_convergence_groups}. Furthermore, using Theorem~\ref{thm:control_concurrent} we can conclude that the steady state behavior of the two clusters is $T$-periodic.

It is straightforward to check that the hypotheses of Theorem \ref{thm:quorum_convergence_groups} are all satisfied if:
\begin{itemize}
\item the biochemical parameters of the two clusters fulfill  the conditions in (\ref{eqn:condition_synchro_max});
\item the coupling function $\phi\left(\cdot\right)$ is increasing.
\end{itemize}

\begin{figure}[thbp]
\begin{center}
\centering \psfrag{x1}[c]{{time (minutes)}}
\centering \psfrag{x2}[c]{{time (minutes)}}
\centering \psfrag{y1}[c]{{$w_i$ (arbitrary units)}}
\centering \psfrag{y2}[c]{{$r\left(t\right)$}}
  \includegraphics[width=8cm]{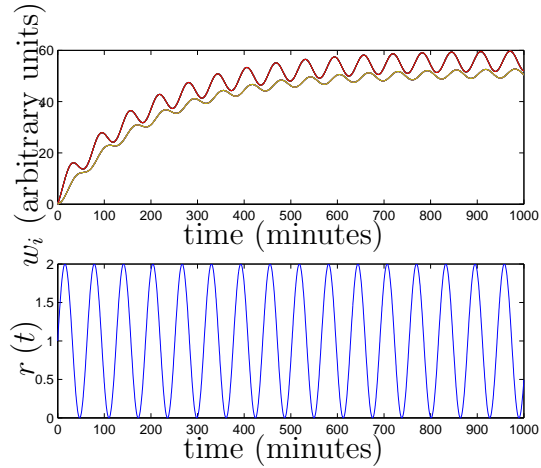}
  \caption{Behavior of (\ref{eqn:math_model_relax_osc_cluster}) when $r\left(t\right)=1+\sin\left(0.1t\right)$. Notice that network nodes have different initial conditions. concurrent synchronization is attained for the network. Both the clusters exhibit a steady state behavior having the same period as $r\left(t\right)$}
  \label{fig: forced_quorum_2}
  \end{center}
\end{figure}

In fact, the topology of the autonomous level of the network is input equivalent by construction. Figure \ref{fig: forced_quorum_2} shows the behavior of (\ref{eqn:math_model_relax_osc_cluster}) when the biochemical parameters of the oscillators are tuned as in the previous Section, and $\phi\left(x\right)= K x$, with $K=0.1$.

\subsubsection{Co-existence of multiple node dynamics}

We analyze the case where the two clusters in the previous Section are
now both connected to a third cluster composed of Van der Pol
oscillators coupled by means of a quorum-sensing mechanism. The three
clusters have three different media, and communication between them
occurs by means of some coupling function. The mathematical model
considered here is then:
\begin{equation}\label{eqn:math_model_relax_osc_cluster_multiple}
\begin{array}{*{20}l}
\dot u_{i1} = \frac{\alpha_1}{1+v_{i1}^\beta} + \frac{\alpha_3 w_{i1}^\eta}{1+w_{i1}^\eta} -d_1 u_{i1}\\
\dot v_{i1} = \frac{\alpha_2}{1+u_{i1}^\gamma} - d_2 v_{i1}\\
\dot w_{i1} = \varepsilon \left(\frac{\alpha_4}{1+u_{i1}^\gamma} - d_3w_{i1}\right)+2 d\left(w_{e1}-w_{i1}\right)\\
\dot w_{e1} = \frac{D_{e}}{N}\sum_{i=1}^N\left(w_{i1}-w_{e1}\right)-d_{e} w_{e1} + \phi\left(w_{e3}\right) - \phi \left(w_{e1}\right)\\
\dot u_{i2} = \frac{\alpha_1}{1+v_{i2}^\beta} + \frac{\alpha_3 w_{i2}^\eta}{1+w_{i2}^\eta} -d_1 u_{i2}\\
\dot v_{i2} = \frac{\alpha_2}{1+u_{i2}^\gamma} - d_2 v_{i2}\\
\dot w_{i2} = \varepsilon \left(\frac{\alpha_4}{1+u_{i2}^\gamma} - d_3w_{i2}\right)+2 d\left(w_{e2}-w_{i2}\right)\\
\dot w_{e2} = \frac{D_{e}}{N}\sum_{i=1}^N\left(w_{i2}-w_{e2}\right)-d_{e} w_{e2}+  \phi\left(w_{e3}\right) - \phi \left(w_{e2}\right)\\
\dot  y_{1i} = y_{2i}\\
\dot y_{2i} = -\alpha\left(y_{1i}^2-\beta\right)y_{2i}-\omega^2y_{1i} +K\left(w_{e3}- y_{1i}\right)\\
\dot w_{e3} = \frac{K}{N_{vdp}}\sum_{i=1}^N\left(y_{2i}-w_{e3}\right)+ g\left(w_{e3}\right)+  \phi\left(w_{e1}\right) + \phi\left(w_{e2}\right) - 2\phi \left(w_{e3}\right)\\
\end{array}
\end{equation}
with $\left[y_{1i}, y_{2i}\right]^T$ denoting the state variables of the $i$-th Van der Pol oscillator, and with $N_{vdp}$ indicating the number of Van der Pol oscillators in the network.
In the above model the Van der Pol oscillators are coupled by means of the medium $w_{e3} \in \R$. The three media, i.e. $w_{e1}$, $w_{e2}$, $w_{e3}$, communicate by means of the coupling function $\phi\left(\cdot\right)$. We assume that the function $g$ governing the intrinsic dynamics of the medium $w_{e3}$ is smooth with bounded derivative. The parameters for the Van der Pol oscillator are set as follows: $\alpha = \beta = \omega = 1$. Notice that now no external inputs is applied on the network.

Recall that Theorem \ref{thm:quorum_convergence_groups} ensures synchronization under the following conditions:
\begin{enumerate}
\item contraction of each cluster composing the network;
\item topology of the autonomous level of the network connected and input equivalent.
\end{enumerate}

Notice that the second condition is satisfied for the network of our interest. Furthermore, contraction of the two clusters composed of genetic oscillators is ensured if the their biochemical parameters satisfy the inequalities in (\ref{eqn:condition_synchro_max}).

To guarantee the convergent behavior of the cluster composed of Van der Pol oscillators, we have to check that there exist two matrix measures, $\mu_\ast$ and $\mu_{\ast\ast}$, showing contraction of the following two matrices:
\begin{subequations} 
\begin{equation}\label{eqn:jac_vdp_1}
J_1=\left[\begin{array}{*{20}c}
0 & 1 \\
-\alpha\left(y_{2i}^2-\beta\right)-\omega^2 & -2\alpha y_{2i}y_{1i}-K
\end{array}\right]
\end{equation}
\begin{equation}\label{eqn:jac_vdp_2}
J_2= \frac{\partial g}{\partial w_{e3}}-K
\end{equation}
\end{subequations}

Now, in \cite{Wan_Slo_05}, by using the Euclidean matrix measure, i.e. $\mu_2$i, it is shown
that the matrix (\ref{eqn:jac_vdp_1}) is contracting if $K>\alpha$. On
the other hand, to ensure contraction of $J_2$, we have to choose
$K>\bar G$, where $\bar G$ is the maximum of $ \frac{\partial
  g}{\partial w_{e3}}$.  Thus, contraction of the cluster composed of
Van der Pol oscillators is guaranteed if the coupling gain, $K$, is
chosen such that:
$$
K > \max\left\{\alpha, \bar G \right\}
$$

\begin{figure}[thbp]
\begin{center}
\centering \psfrag{x1}[c]{{time (minutes)}}
\centering \psfrag{x2}[c]{{time (minutes)}}
\centering \psfrag{y1}[c]{{$w_{i1}$, $w_{i2}$}}
\centering \psfrag{y2}[c]{{$y_{2i}$}}
  \includegraphics[width=8cm]{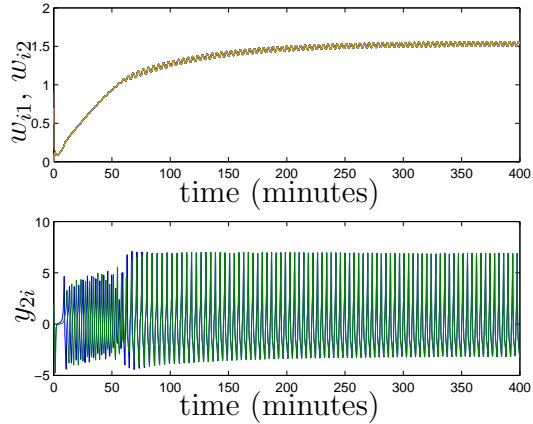}
  \caption{Behavior of (\ref{eqn:math_model_relax_osc_cluster_multiple}). The two clusters of genetic oscillators (not directly connected) converge onto the same common evolution. Synchronization for the cluster composed of Van der Pol oscillators is also attained.}
  \label{fig:cluster_vdp}
  \end{center}
\end{figure}

In Figure \ref{fig:cluster_vdp}, we set $g\left(x\right) =\sin\left(x\right)$, $K=2.5$, $N_{vdp}=2$ and $\phi\left(x\right)= Kx$, with $K=3$. Such a Figure shows that concurrent synchronization of (\ref{eqn:math_model_relax_osc_cluster_multiple}) is attained, in agreement with the theoretical analysis.

\subsection{Analysis of a general Quorum-Sensing pathway}

In the previous Section, we showed that our results (with appropriate
choice of matrix measure) can be used to derive easily verifiable
conditions on the biochemical parameters of the genetic oscillator
ensuring contraction, and hence synchronization (onto a periodic orbit
of desired period) and concurrent synchronization. We now show that
our methodology can be applied to analyze a wide class of biochemical
systems involved in cell-to-cell communication.

\begin{figure}[thbp]
\begin{center}
  \includegraphics[width=10cm]{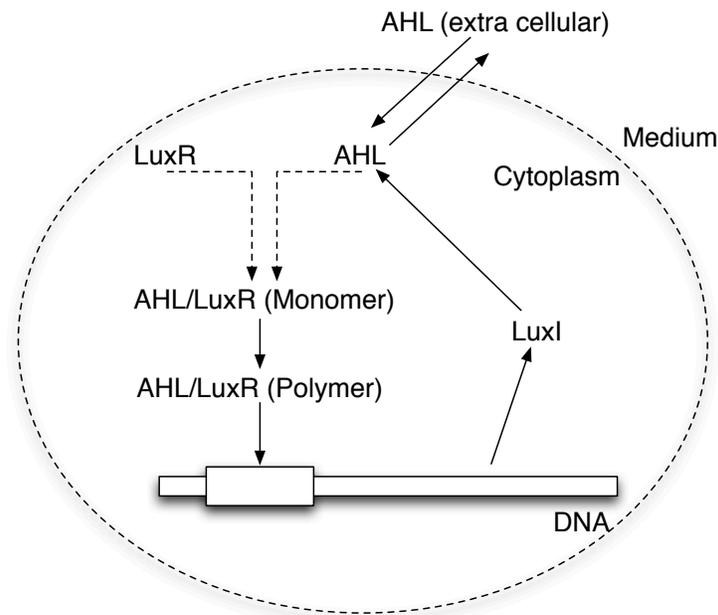}
  \caption{The quorum sensing pathway implemented by $AHL$}
  \label{fig:pathway_quorum}
  \end{center}
\end{figure}

We focus on the analysis of the pathway of the quorum sensing mechanism that uses as autoinducers, molecules from the AHL (acyl homoserine lactone) family. The quorum sensing pathway implemented by AHL (see Figure \ref{fig:pathway_quorum}) is one of the most common for bacteria and drives many transcriptional systems regulating their basic activities.

We now briefly describe the pathway of our interest (see \cite{Mul_Kut_Hen_Rot_Har_06} for further details). The enzyme LuxI produces AHL at (approximately) a constant rate. AHL in turn diffuses into and out of the cell and forms (in the cytoplasm) a complex with the receptor LuxR. Such complex polymerizes and then acts as a transcription factor, by binding the DNA. This causes the increase of the production of LuxI, generating a positive feedback loop.

The pathway can be described by a set of ordinary differential equations (using the law of mass action, see \cite{Doc_Kee_04}, \cite{Mul_Kut_Hen_Rot_Har_06}). Specifically, denoting with $x_e$ the mass of AHL outside of the cell and with $x_c$ the mass of AHL within the cell, we have the following mathematical model:
\begin{equation}\label{eqn:quorum_pathway}
\begin{array}{*{20}l}
\dot x_c = \alpha + \frac{\beta x_c^n}{x_{thresh}^n + x_c^n}- \gamma_cx_c -  d_1 x_c -  d_2 x_e\\
\dot x_e =  d_1 x_c -  d_2 x_e -\gamma_e x_e\\
\end{array}
\end{equation}
The physical meaning of the parameters in (\ref{eqn:quorum_pathway}) is given in Table \ref{tab:parameters_pathway}.

\begin{table}[th] 
\caption{Biochemical parameters for system (\ref{eqn:quorum_pathway})}
\centering
\label{tab:parameters_pathway}
\begin{tabular}{|c| c|} 
\hline
Parameter & physical meaning\\
\hline
$\alpha$ & Low production rate of $AHL$\\
$\beta$ & Increase of production rate of $AHL$\\
$\gamma_c$ & Degradation rate of $AHL$ in the cytosol\\
$\gamma_e$ & Degradation rate of $AHL$ outside the cell\\
$ d_1$ & Diffusion rate of the extracellular $AHL$\\
$ d_2$ & Diffusion of the intracellular $AHL$\\
$x_{thresh}$ & Threshold of $AHL$ between low and increased activity\\
$n$ & Degree of polymerization\\
\hline
\end{tabular}
\end{table}

Now, contraction of the above system is guaranteed if 
\begin{enumerate}
\item $-\gamma_c + \frac{2\beta x_{thresh}^2x_c}{\left(x_{thresh}^2+x_c^2\right)^2}$ is uniformly negative definite;
\item $- d_2 - \gamma_e$ is uniformly negative definite.
\end{enumerate}
Recall that $x_c$ and $x_e$ are both scalars.
Now, the second condition is satisfied since system parameters are all positive. That is, to prove contraction we have only to guarantee that 
$$
-\gamma_c + \frac{2\beta x_{thresh}^2x_c}{\left(x_{thresh}^2+x_c^2\right)^2}
$$
is uniformly negative. Since
$$
-\gamma_c + \frac{2\beta x_{thresh}^2x_c}{\left(x_{thresh}^2+x_c^2\right)^2}\le -\gamma_c + \frac{3\beta \sqrt{3}}{8 x_{thresh}}
$$
contraction is ensured if the biochemical parameters $\beta$, $g$ and $x_{thresh}$ fulfill the following condition
$$
\frac{\beta}{x_{thresh}} <\frac{8\gamma_c}{3\sqrt{3}}
$$

\section{Concluding remarks}

In this paper, we presented a systematic methodology to derive
conditions for the global exponential convergence of biochemical
models modeling quorum sensing systems. To illustrate the
effectiveness of our results and to emphasize the use of our
techniques in synthetic biology design, we analyzed a set of
biochemical networks where the quorum sensing mechanism is involved as
well as a typical pathway of the quorum sensing. In all such cases we
showed that our results can be used to determine system parameters and
dynamics ensuring convergence.

\appendix

\section{Proofs}

To prove Theorem \ref{thm:network_convergence_multiple_nodes} we need the following Lemma, which is a generalization of a result proven in \cite{Rus_diB_Son_09}:

\begin{Lemma}\label{lem:partition}
Consider the block- partition for a square matrix $J$:
$$
J =\left[\begin{array}{*{20}c}
A(x) & B(x,y)\\
C(x,y) & D(y)\\
\end{array}\right]
$$
where $A$ and $D$ are square matrices of dimensions $n_A \times n_A$ and $n_D\times n_D$ respectively. Assume that $A$ and $B$ are contracting with respect to $\mu_A$ and $\mu_D$ (induced by the vector norm $\abs{\bullet}_A$ and $\abs{\bullet}_D$). Then, $J$ is contracting if there exists two positive real numbers $\theta_1$, $\theta_2$ such that
$$
\begin{array}{*{20}c}
\mu_A(A)+\frac{\theta_2}{\theta_1}\norm{C(x,y)}_{A,D} \le - c_A^2\\
\mu_D(D)+\frac{\theta_1}{\theta_2}\norm{B(x,y)}_{D,A} \le - c_B^2\\
\end{array}
$$
where $\norm{\bullet}_{A,D}$ and $\norm{\bullet}_{D,A}$ are the operator norms induced by $\abs{\bullet}_A$ and $\abs{\bullet}_D$ on the linear operators $C$ and $B$.  Furthermore, the contraction rate is $c^2 =\max \left\{c_A^2, c_B^2\right\}$.
\end{Lemma}
\proof
Let $z:=(x,y)^T$. We will show that, with the above hypotheses, $J$ is contracting with respect to the matrix measure induced by the following vector norm:
$$
\abs{z} := \theta_1 \abs{x}_A + \theta_2\abs{y}_D
$$
with $\theta_1, \theta_2 >0$. In this norm, we have
$$
\abs{(I+hJ)z} = \theta_1\abs{(I+hA)x + hBy}_A + \theta_2\abs{(I+hD)y + hCx}_D
$$
Thus,
$$
 \abs{(I+hJ)z} \le \theta_1\abs{(I+hA)x}_A + h\theta_1\abs{By}_{D,A} + \theta_2\abs{(I+hD)y}_D + h\theta_2\abs{Cx}_{A,D}
$$
Pick now $h>0$ and a unit vector $z$ (depending on $h$) such that $\norm{(I+hJ)z}= \abs{(I+hJ)z}$. We have, dropping the subscripts for the norms:
$$
\frac{1}{h}(\norm{I+hJ}-1 ) \le \frac{1}{h}\left(\norm{I+hA}-1+\frac{\theta_2}{\theta_1}h\norm{C}\right)\abs{x}\theta_1 +  \frac{1}{h}\left(\norm{I+hD}-1+\frac{\theta_1}{\theta_2}h\norm{B}\right)\abs{y}\theta_2
$$
Since $1=\abs{z} = \theta_1 \abs{x}_A + \theta_2 \abs{y}_B$, we finally have
$$
\frac{1}{h}(\norm{I+hJ}-1 ) \le \max \left\{ \frac{1}{h}\left(\norm{I+hA}-1+\frac{\theta_2}{\theta_1}h\norm{C}\right),   \frac{1}{h}\left(\norm{I+hD}-1+\frac{\theta_1}{\theta_2}h\norm{B}\right)\right\} 
$$
Taking now the limit for $h\rightarrow 0^+$:
$$
\mu\left(J\right) \le \max\left\{ \mu(A)_A +\frac{\theta_2}{\theta_1}\norm{C}, \mu(D)_D +\frac{\theta_1}{\theta_2}\norm{B} \right\}
$$
thus proving the result.
\endproof

Following the same arguments, Lemma \ref{lem:partition} can be straightforwardly extended to the case of a real matrix $J$ partitioned as
$$
J =\left[\begin{array}{*{20}c}
J_{11} & J_{12} & \ldots & J_{1N}\\
\ldots & \ldots & \ldots & \ldots \\
J_{N1} & J_{N2} & \ldots & J_{NN}\\
\end{array}\right]
$$
where the diagonal blocks of $J$ are all square matrices. Then  $J$ is contracting if
\begin{equation}\label{eqn:lemma_useful}
\begin{array}{*{20}l}
\mu(J_{11})+\frac{\theta_2}{\theta_1}\norm{J_{12}}+ \ldots + \frac{\theta_N}{\theta_1}\norm{J_{1N}} \le - c_{11}^2\\
\ldots \\
\mu(J_{NN})+\frac{\theta_1}{\theta_N}\norm{J_{N1}}+ \ldots + \frac{\theta_{N-1}}{\theta_N}\norm{J_{1N}} \le - c_{NN}^2\\
\end{array}
\end{equation}
(where subscripts for matrix measures and norms have been neglected).

\subsection*{Proof of Theorem \ref{thm:network_convergence_multiple_nodes}}

The assumption of input equivalence for the nodes implies the
existence of a linear invariant subspace associated to the concurrent synchronization steady state regime. We
will prove convergence towards such a subspace, by proving that the
network dynamics is contracting. Let $\mu_f$ be the matrix measure
where the nodes dynamics is contracting and define: $X:=(x_1^T,\ldots,
x_N^T)^T$, $F(X)$ as the stack of all intrinsic nodes dynamics, $H(X)$ the stack
of nodes coupling functions. We want to prove that there exist a
matrix measure, $\mu$, (which is in general different from $\mu_f$)
where the whole network dynamics is contracting. Denote with $L:=\left\{l_{ij}\right\}$ the Laplacian matrix \cite{God_Roy_01} of the network and define the matrix $\tilde L(X)$, whose $ij$-th block, $\tilde L_{ij}(X)$, is defined as follows:
$$
\tilde L_{ij}(X) := l_{ij} \frac{\partial h_{\gamma(i)}}{\partial x_j}
$$
(Notice that if all the nodes are identical and have the same dynamics and the same coupling functions, then $\tilde L$ can be written in terms of the Kronecker product, $\otimes$, as $(L\otimes I_n)\frac{\partial H}{\partial X}$, with $n$ denoting the dimension of the nodes and $I_n$ the $n\times n$ identity matrix.)

 The Jacobian of
(\ref{eqn:net_multiple_dynamics}) is then:
\begin{equation}\label{eqn:net_jacobian}
J:=\left[\frac{\partial F}{\partial X} - \tilde L(X)\right]
\end{equation}
The system is contracting if
$$
\mu\left(\frac{\partial F}{\partial X} -\tilde L(X) \right)
$$
is uniformly negative definite.
Now:
$$
\mu\left(\frac{\partial F}{\partial X} - \tilde L(X)  \right)
 \le \mu\left(\frac{\partial F}{\partial X}\right) + \mu\left(-\tilde L(X) \right)
$$
Notice that, by hypotheses, the matrix $- \tilde L(X)$ has negative diagonal blocks and zero column sum. Thus, using (\ref{eqn:lemma_useful}) with $\theta_i=\theta_j$ for all $i,j=1,\ldots, N$, $i \ne j$ yields
$$
\mu\left(- \tilde L(X)\right) = 0
$$
Thus:
$$
\mu\left(\frac{\partial F}{\partial X} - \tilde L(X)  \right)
 \le \mu\left(\frac{\partial F}{\partial X}\right)
$$
Since the matrix $\frac{\partial F}{\partial X}$ is block diagonal, i.e. all of its off-diagonal elements are zero, (\ref{eqn:lemma_useful}) yields:
$$
 \mu\left(\frac{\partial F}{\partial X}\right) = \max_{x,t,i}\left\{\mu_f\left(\frac{\partial f_{\gamma(i)}}{\partial x}\right)\right\}
$$
The theorem is then proved by noticing that by hypothesis the right
hand side of the above expression is uniformly negative.


 \end{document}